%% file: main.tex
\documentclass{article} 
\usepackage{graphicx} 
\usepackage{amsmath}
\usepackage{amsfonts}
\usepackage{cite}
\usepackage{hyperref}
\usepackage{diagbox}
\usepackage{multirow}
\usepackage{adjustbox}
\usepackage{longtable}
\usepackage{comment}
\usepackage{xcolor}
\usepackage[left=2cm, right=2cm, top=3cm]{geometry}
\usepackage{amsthm}
\usepackage{authblk}

\usepackage{subcaption} 

\DeclareMathOperator{\diag}{diag}
\DeclareMathOperator{\rank}{rank}
\DeclareMathOperator{\mrank}{mrank}
\DeclareMathOperator{\mem}{mem}
\DeclareMathOperator{\conv}{conv}

\AtBeginDocument{}

\title{Mosaic-skeleton approximation is all you need for Smoluchowski equations}
\author[1,2]{Roman Dyachenko}
\author[2,3]{Sergey Matveev}
\author[2,3]{Bulat Valiakhmetov}
\affil[1]{HSE University and Skolkovo Institute of Science and Technology, Moscow, Russia}
\affil[2]{Marchuk Institute of Numerical Mathematics, RAS, Moscow, Russia}
\affil[3]{Lomonosov MSU, faculty of Computational Mathematics and Cybernetics, Moscow, Russia}

\begin{document}


\maketitle
\begin{abstract}
In this work we demonstrate a surprising way of exploitation of the mosaic--skeleton approximations for efficient numerical solving of aggregation equations with many applied kinetic kernels. The complexity of the evaluation of the right-hand side with $M$ nonlinear differential equations basing on the use of the mosaic-skeleton approximations is $\mathcal{O}(M \log^2 M)$ operations instead of $\mathcal{O}(M^2)$ for the straightforward computation. The class of kernels allowing to make fast and accurate computations via our approach is wider than analogous set of kinetic coefficients for effective calculations with previously developed algorithms. This class covers the aggregation problems arising in modelling of sedimentation, supersonic effects, turbulent flows, etc. We show that our approach makes it possible to study the systems with $M=2^{20}$ nonlinear equations within a modest computing time.
\end{abstract}

\input{intro.tex}
\input{approximation_kernels.tex}

\input{calculation_operators.tex}
\input{numerical_experiments.tex}

\section{Conclusion and future work}

In this work we demonstrate that the mosaic-skeleton matrix format allows to approximate the broad family of coagulation kernels. 
The novel method allows to evaluate the Smoluchowski operator numerically within $\mathcal{O}(RM \log^2 M)$ operations.
It make possible to simulate aggregation process even for the wide particle size distributions and for the kernels with formally high ranks.
We incept these new fast methods into the classical Runge-Kutta methods and verify their accuracy in comparison with the Monte Carlo \texttt{LRMC} approach.
The tests of our approach are presented with using up to $M = 2^{20}$ nonlinear equations.

The developed methodology is implemented as an open-source software and allows the researchers to use our approach as a blackbox solver, setting only the coagulation kernel and initial conditions. 
In the future, we plan to study if the mosaic-skeleton approximations could be also useful for accelerating the Monte Carlo methods.

\section*{Acknowledgements}

This work is supported by Russian Science Foundation, Project 25--21--00047 (the project webpage is available \url{https://rscf.ru/project/25-21-00047/}). Authors are grateful to Alexander Osinsky for useful discussions.

\medskip

\bibliographystyle{ieeetr}
\bibliography{referencesMSCSmol}

\end{document}

%% file: intro.tex
\section{Introduction}
%
Aggregation is a basic natural process and plays an important role in numerous phenomena \cite{krapivsky2010kinetic} including the processes at the microscale (e.g., polymerization \cite{Flory41, skorych2019investigation} or dust formation \cite{friedlander2000smoke}) and macroscale (e.g., particle formation at Saturn's rings \cite{matveev2017oscillations}). Series of useful facts about growing random graphs can be found in the framework of aggregation equations \cite{ER60, Janson93, Lushnikov05, aldous1999deterministic}. All in all, this knowledge might be useful for revealing the mechanisms of the community formation in the social networks and biological systems \cite{dorogovtsev2003evolution}. 

We study the aggregation process assuming that it is a collision-controlled process and that particles fill the space homogeneously. In such case, one may pay interest not to evolution of the distinct particles but to concentrations $n_s(t)$ of the cluster of size $s$ per unit volume of the system. In case of pairwise collisions leading to the mass-conserving aggregation with rates $K_{ij} = K_{ji}\geq 0$ (also called kernels) the celebrated Smoluchowski equations can be derived \cite{Smol17, smoluchowski1916drei}:
\begin{equation} \label{eq:smoluch_euler}
    \frac{d n_s}{dt} = \frac{1}{2} \sum\limits_{i+j=s} K_{ij} n_i n_j - \sum\limits_{j = 1}^{\infty} K_{sj} n_s n_j.
\end{equation}
These balance equations correspond to the irreversible aggregation process and can be studied analytically only for very special cases of kernels and initial conditions \cite{Leyvraz03}. In most applications, researchers utilize the numerical methods that allow one to approximate the solution with some finite accuracy. Today, there exist many fundamentally different approaches that allow to study the aggregation equations with the use of computers \cite{singh2022challenges}. Namely, they can be finite-difference \cite{matveev2015fast, osinsky2020low, hackbusch2006efficient}, finite-volume \cite{singh2021accurate}, direct simulation Monte Carlo approaches \cite{sabelfeld2018hybrid, kruis2000direct, eibeck2000efficient, patterson2011stochastic, goodson2002efficient}, the specific coarse graining tricks \cite{lee2000validity, stadnichuk2015smoluchowski} and even the homotopy perturbation methods \cite{yadav2023note}. Each of these methods has its specific advantages and drawbacks. For example, the direct simulation Monte Carlo methods naturally fulfill the mass conservation laws, but their accuracy is limited. On the other hand, the deterministic finite-difference approach is very accurate but might be resource-consuming. Coupling of the finite-difference methods with the ideology of the low-rank decomposition allows to partially solve this problem and perform extremely accurate simulation within modest times \cite{matveev2015fast, osinsky2020low, hackbusch2006efficient}. However, there are still classes of aggregation kernels that do not have exactly or numerically low rank. Some examples of these kernels are the following:
\begin{align}
    K_{i j} &= \left\lbrace \begin{matrix} \left(i^{1/3} +j^{1/3} \right)^2 |i^{2/3}-j^{2/3}|, & i \neq j, \\
    ~ & ~\\
    \quad (i^{1/3} + j^{1/3})(i^{-1/3} + j^{-1/3}), & i = j. \end{matrix} \right. \label{eq:kernel_stream} \\
    K_{i j} &= \left\lbrace \begin{matrix} \dfrac{(i+j)\left(i^{1 / 3}+j^{1 / 3}\right)^{2 / 3}}{(i j)^{5 / 9}\left|i^{2 / 3}-j^{2 / 3}\right|}, & i \neq j, \\
    ~ & ~\\
    \quad (i^{1/3} + j^{1/3})(i^{-1/3} + j^{-1/3}), & i = j. \end{matrix} \right. \label{eq:kernel_baikal}
\end{align}
Such kernels (see e.g. \eqref{eq:kernel_stream}) may arise in problems with spatially inhomogeneous aggregation of particles moving within some stream. Surprisingly, the popular and fine-tuned Monte Carlo approaches also degenerate in terms of efficiency for the problems with such kernels. It motivates us to revisit and reformulate the methodology utilizing low-rank matrix structures for the more general class of matrices with low mosaic rank.

In our work, we extend the previously proposed approach \cite{koch2008h} through the application of the adaptive cross methods \cite{goreinov2001maximal, goreinov2010find, zheltkov2015parallel, bebendorf2000approximation} for the mosaic approximation of the kinetic coefficients. This idea allows us to solve huge systems of up to $2^{20}$ equations with no need for supercomputing facilities. We also confirm the efficiency of our approximations for a large number of kernels encountered in practice and conduct the convergence analysis in comparison with the family of Monte Carlo methods.

The main contributions of our work are:
\begin{itemize}
    \item great compression rate on a broad family of coagulation kernels within the mosaic--skeleton format;
    \item new algorithms for efficient evaluation of the Smoluchowski operator based on the mosaic-skeleton structure;
    \item incorporation of the adaptive cross approximation allowing to deal with huge systems of aggregation equations;
    \item an open high-performance solver for solving Smoluchowski equations \footnote{\href{https://github.com/DrEternity/FDMSk-Smoluchowski}{https://github.com/DrEternity/FDMSk-Smoluchowski}}.
\end{itemize}

%% file: approximation_kernels.tex
\section{Approximation of kernels} \label{sec:Approximation of kernels}
Let us describe the scheme for modeling the coagulation process using (\ref{eq:smoluch_euler}). We first consider a coagulation equation with the particle size limit set as $M$:
\begin{equation} \label{eq:smoluch_euler_approximation}
    \frac{dn_s}{dt} = \underbrace{\frac12 \sum\limits_{i+j=s} K_{ij} n_i n_j}_{f_1(s)} - \underbrace{\sum\limits_{j = 1}^{M} K_{sj} n_s n_j}_{f_2(s)}, \quad s = 1, \ldots,  M.
\end{equation}
Such system of ODEs can be solved with various time-integration approaches.
Nevertheless, in order to perform the time-step one has to calculate the operators $f_1(s)$, $f_2(s)$ for all $s =1, 2, \ldots, M$.

It is easy to see that the direct calculation requires $\mathcal{O}(M^2)$ arithmetic operations for both $f_1$ and $f_2$ terms.
However, in recent years, some efficient methods have been developed~\cite{matveev2015fast, osinsky2020low, chaudhury2014computationally}.
These approaches utilize low--rank approximations of the kernel to speed up the calculation of $f_1$ and $f_2$.


\subsection{Low-rank method} \label{sec:lowrank_kernels}

In this section we show a certain approach allowing to speedup computation of the operators  $f_1$ and $f_2$ from Eq.~\eqref{eq:smoluch_euler_approximation}. At first, we construct a low-rank decomposition of the coagulation kernel:
\begin{align}
    K_{ij} &\approx \sum_{\alpha=1}^R U_{i \alpha} V_{j \alpha}, \label{eq:kernel_split} \\
    K &\approx U V^\top, \label{eq:skelet_decomp}
\end{align}
where $U, V \in \mathbb{R}^{M \times R}$. Any relevant approach can be used for construction of $U$ and $V$ e.g. the singular value decomposition (SVD) \cite{horn2012matrix}, its randomized version \cite{halko2011finding} or the adaptive cross approximation \cite{goreinov2001maximal, goreinov2010find, zheltkov2015parallel, bebendorf2000approximation}. 

The adaptive cross method allows one to deal with function-generated matrices and has a linear complexity $\mathcal{O}(MR^2)$ with respect to the matrix size $M$. These two features allow us to approximate kernels for systems with up to $M = 2^{20}$ equations.


As soon as the low-rank decomposition is constructed, we obtain the relation for the first operator:
\begin{equation*}
    f_1(s) = \frac12\sum\limits_{i+j=s} K_{ij} n_i n_j = \frac12\sum \limits_{i+j=s} \sum_{\alpha=1}^R U_{i\alpha} V_{j\alpha} n_i n_j = \frac12\sum_{\alpha=1}^R \sum \limits_{i+j=s} \left[ U_{i\alpha} n_i \right] \left[V_{j\alpha} n_j \right].
\end{equation*}
It is the sum of $R$ vector convolutions:
\begin{equation*}
    f(s) = \sum\limits_{i+j=s} g(i) h(j).
\end{equation*}
Each summand can be calculated within $\mathcal{O}(M \log M )$ operations using the fast Fourier transform. 

For the second operator, we obtain the following relation:
\begin{equation*}
    f_2(s) = n_s\sum\limits_{j = 1}^{M} K_{sj} n_j = n_s \sum\limits_{j = 1}^{M} \sum_{\alpha=1}^R U_{s\alpha} V_{j\alpha} n_j = n_s \sum_{\alpha=1}^R U_{s\alpha} \sum\limits_{j = 1}^{M} V_{j\alpha} n_j.
\end{equation*}
Or in matrix form it can be written like this:
\begin{equation*}
    f_2 = n \odot \left( U V^\top n\right),
\end{equation*}
where $U$ and $V$ are the factors of the skeleton decomposition, $n$ is an input vector and $\odot$ is an element-wise product. Evaluation of the whole vector $f_2$ requires $\mathcal{O}(MR)$ operations.

Thus, we obtain a numerical approach for calculation of the Smoluchowski operator within $\mathcal{O}(R M\log M )$ operations instead of the initial $\mathcal{O}(M^2)$. 

Final structure of the numerical method exploiting the low-rank representations of the aggregation coefficients requires some time-integration scheme. In this work, we utilize the classical Runge--Kutta--Fehlberg method~\cite{fehlberg1970classical, matveev2024adaptive}.

The quality of the low-rank approximations of matrix $K$ in the Frobenius norm corresponds to the decrease of its singular values.  Let $A = U\Sigma V^{\top} \in \mathbb{R}^{m \times n}$, $m \leq n$, be the singular value decomposition of $A$, where $\Sigma=\diag(\sigma_1,\ldots,\sigma_m)$ and\; $U, V^T$ -- orthogonal matrices. Taking the truncated matrices $U_r$, $V_r$ as the first $r$ columns of $U$ and $V$ and truncating $\Sigma_r = \diag(\sigma_1,\ldots,\sigma_r)$, we obtain the rank-$r$ matrix $A_r = U_r\Sigma_r V_r^{\top}$. $A_r$ is an optimal approximation of $A$ in terms of the Frobenius norm:
    $$\| A - A_r\|_F = \min_{\operatorname{rank}(B) \leq r} \|A-B\|_{\text{F}} = \sqrt{\sigma^2_{r+1} + \cdots + \sigma^2_m}.$$
This expression for the value of the approximation error is known as the famous Eckart-Young-Mirsky theorem (see e.g. ~\cite{horn2012matrix}).

Hence, if the ``tail'' $\sigma_{r+1}$, $\ldots$, $\sigma_m$ of the singular values decreases rapidly, then the whole matrix can be accurately approximated with a low-rank one. The family of such matrices is rather broad and includes a lot of function-generated matrices \cite{goreinov1997theory} and also positive semidefinite Hankel matrices \cite{beckermann2017singular}.

An approach basing on low-rank matrix approximations can be applied successfully to a wide class of applied phenomena related to aggregation kinetics. The list of possible kernels includes the ballistic~\cite{carnevale1990statistics, trizac2003correlations, paul2018dimension}, Brownian~\cite{lai1972self, thorn1994dynamic, odriozola2004irreversible}, and others \cite{dyachenko2023finite, matveev2017oscillations, skorych2019investigation}. However, many kernels (for example given by Eqs.~\eqref{eq:kernel_stream}~and~\eqref{eq:kernel_baikal}) ~\cite{liao2010literature, sheng2007simulation, floyd2014soot, derevich2007coagulation, osinsky2024hydrodynamic, kostoglou2020critical} cannot be approximated via the low--rank format.

\begin{figure}[ht]
     \centering
     \begin{subfigure}[t]{0.49\textwidth}
         \centering
         \includegraphics[width=\textwidth]{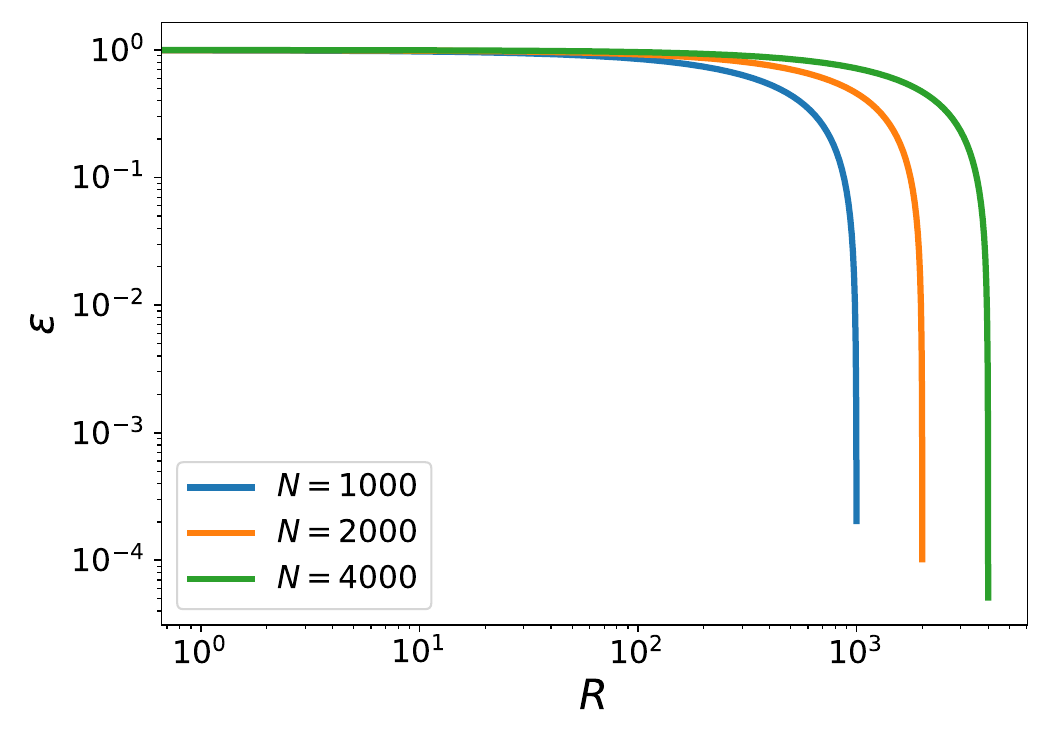}
            \caption{$K_{ij}=\dfrac{(i+j)\left(i^{1/3}+j^{1/3}\right)^{2 / 3}}{( i j)^{5/9}\left|i^{2/3}-j^{2/3}\right|}$}
     \end{subfigure}
     \hfill
     \begin{subfigure}[t]{0.49\textwidth}
         \centering
         \includegraphics[width=\textwidth]{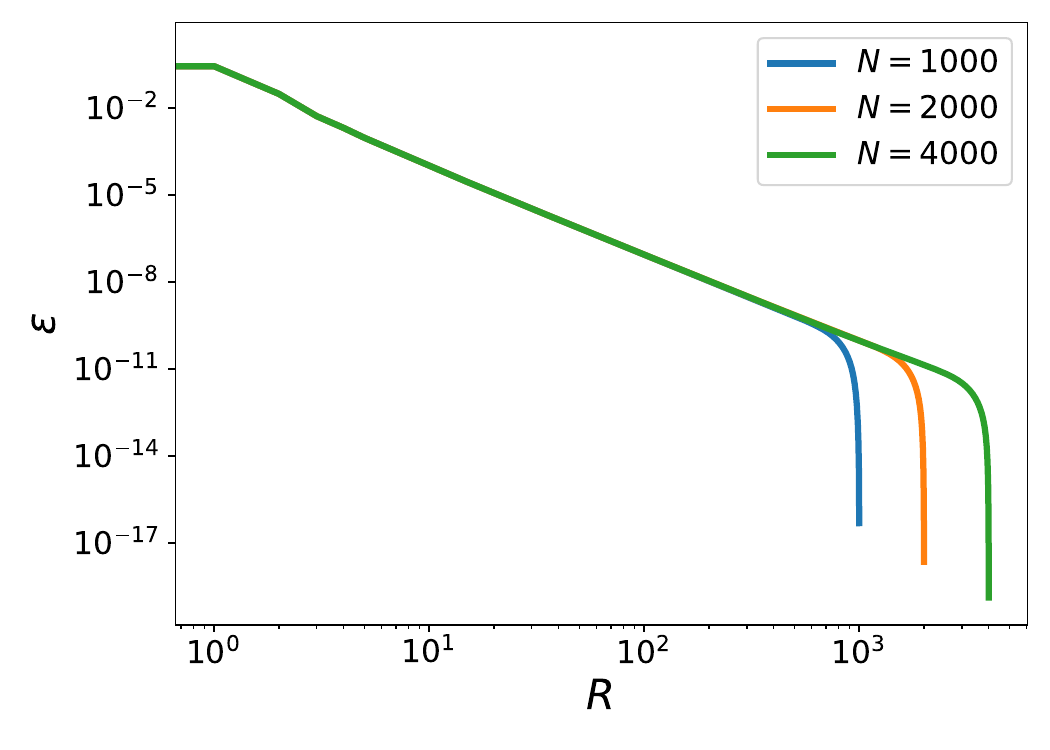}
         \vspace{-3pt}
         \caption{$K_{ij} = \left(i^{1/3} +j^{1/3} \right)^2 |i^{2/3}-j^{2/3}|$}
     \end{subfigure}
        \caption{Smoluchowski kernels low--rank approximation error in the Frobenius norm depending on the rank $R$ for different matrix sizes $N$}
        \label{fig:bad_sing_values}
\end{figure}

The slow decrease of the singular values of the kernels \eqref{eq:kernel_stream} and \eqref{eq:kernel_baikal} can be seen on Figure~\ref{fig:bad_sing_values}. This trouble inspires us to extend the existing approach  using the wider class of low parametric matrices for solving the Smoluchowski-class equations.


\subsection{Mosaic--skeleton approximations} \label{sec:mosaic_skeleton_kernels}

Mosaic--skeleton matrix representation emerged as a way allowing to approximate function generated matrices, in particular the matrices related to the integral operators~\cite{voevodin1979obodnom, hackbusch1989fast}.
This format is also known as a hierarchical format~\cite{hackbusch2015hierarchical} or $\mathcal{H}$-matrix representation.

In general, the matrix $K$ can be considered as the table of some function values:
\begin{equation*}
    K_{ij} = f(y_i, x_j).
\end{equation*}
If $f(y, x)$ is smooth enough and two sets of values $Y = \{y_i\}$, $X = \{x_j\}$ correspond to some quasi-uniform spatial grids, then $K$ may be well approximated within the low--rank format~\eqref{eq:skelet_decomp}. 
In particular, it is known that for the \textit{asymptotically smooth} functions the corresponding matrices can also be described with modest number of parameters growing almost linearly with respect to its size~\cite{tyrtyshnikov1996mosaicskeleton}. Hence, the mosaic--skeleton format may be considered as an extension of basic low-rank representation.

\begin{figure}[ht]
     \centering
     \begin{subfigure}[t]{0.47\textwidth}
         \centering
         \includegraphics[width=\textwidth]{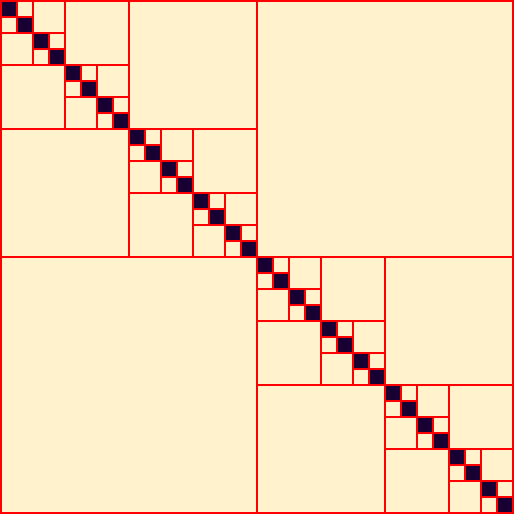}
         \caption{Dense diagonal}
         \label{fig:different_splitting:monodiag}
     \end{subfigure}
     \hfill
     \begin{subfigure}[t]{0.47\textwidth}
         \centering
         \includegraphics[width=\textwidth]{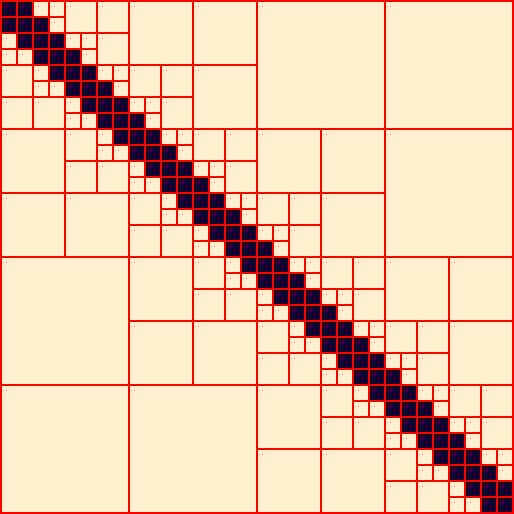}
         \caption{Dense diagonal, sub-diagonal and super-diagonal}
          \label{fig:different_splitting:tridiag}
     \end{subfigure}
        \caption{Mosaic partitioning depending on the low--rank criteria}
        \label{fig:different_splitting}
\end{figure}

Let us consider the whole matrix describing the full set of pairwise interactions between entries of sets $Y$ and $X$ defined by the elements of $K$.
For the sake of finding the low-rank structure, we can split each of them into two smaller subsets and check its presence at four new blocks.
If some of them cannot be approximated, then we split those ones recursively.
One has to continue this splitting procedure until each new block can be approximated precisely via basic skeleton format or its size is small enough.
These small blocks have to be evaluated and stored fully, as dense matrices.
Finally, the matrix~$K$ is split into blocks of different sizes and formats (either dense or the low-rank ~\eqref{eq:skelet_decomp}).
We name such blockwise structure the mosaic--skeleton format and present it schematically in Figure~\ref{fig:different_splitting}.

There are some key questions to clarify such procedure:
\begin{enumerate}
\item
    How to split blocks (and corresponding spatial points of sets $Y$ and $X$)?    
\item
    How to determine whether the block has numerically low rank or not?
\item
    How to get the approximation factors $U$ and $V$ for every block of such kind?
\end{enumerate}

Complete implementation details and nuances of the general algorithm could be found in different papers~\cite{valiakhmetov2023msk-ruscd, tyrtyshnikov2000incomplete, hackbusch2015hierarchical}.
In our case, we utilize the basic version of the algorithm constructing the mosaic-skeleton approximation.

Since the aggregation kernel elements $K$ depend on indices in a straight-forward way, the points of $Y$ and $X$ are set as $x_i \equiv y_i \equiv i \in \{1, 2, \ldots, M\}$.
The questions listed above get the following answers:
\begin{enumerate}
\item
    For each block its splitting (if necessary) is provided just by halves of rows and columns.
    It produces four quarter blocks.
\item
    We choose the following criteria allowing to determine the low-rank blocks: (a) only diagonal blocks are marked as dense; (b) only diagonal, sub-diagonal and super-diagonal blocks are dense.
    All blocks except these are treated as low-rank. 
    These two options lead to different mosaic partitionings depicted on Figure~\ref{fig:different_splitting}.
\item
    We construct the approximation of blocks that are expected to have numerically low rank using the adaptive cross approximation procedure described in the paper~\cite{zheltkov2015parallel}. 
    This algorithm allows to avoid the storage of full submatrices in memory, and evaluates the approximation ``on-the-fly''.
\end{enumerate}

In Table~\ref{table:mosaic_ranks} we present some examples of kernels from different studies with their compression rate and maximal rank of a block in the mosaic-skeleton format.
In these experiments we consider both mosaic partitioning options.
Approximation is obtained with relative accuracy $\varepsilon = 10^{-6} \text{ and } 10^{-12}$ in the Frobenius norm. 
The final compression rate is measured with respect to dense matrix with $M^2$ entries.

The described approach for approximation of the aggregation kernels demonstrates an excellent result. In the worst case, the rank $8$ is enough to achieve a relative accuracy of $10^{-12}$, for each complex kernel from Table~\ref{table:mosaic_ranks}. These results show the versatility of our methodology and allow to apply our approach as a black-box solution for numerical modelling with a wide class of the aggregation kernels.

\begin{table}[ht!]
\begin{tabular}[c]{| c || c | c | c | c | p{3cm} |}
     \hline
    \multirow{2}{*}{\diagbox{$K_{ij}$}{$\varepsilon$}} & \multicolumn{2}{c|}{Partitioning~1} & \multicolumn{2}{c|}{Partitioning~2} & \multirow{2}{*}{Source} \\
    & $10^{-6}$ & $10^{-12}$ & $10^{-6}$ & $10^{-12}$ &  \\    \hline
    \hline
    \rule{0pt}{1.5em}$\left(i^{1/3} +j^{1/3} \right)^2 |i^{2/3}-j^{2/3}|$& $5 (0.01)$ & $5 (0.01)$ & $4 (0.03)$& $5 (0.06)$ & aggregation within flow \cite{horvai2008coalescence, casamatta1985modelling, friedlander2000smoke, zagidullin2022aggregation}\\
    \hline
    \rule{0pt}{1.5em}$\dfrac{(i+j)\left(i^{1 / 3}+j^{1 / 3}\right)^{2 / 3}}{(i j)^{5 / 9}\left|i^{2 / 3}-j^{2 / 3}\right|}$ & $15 (0.04)$ & $34 (0.07)$ & $4 (0.04)$& $8 (0.06)$ & baseline example\\
    \hline
    \rule{0pt}{1.5em} $(i^{1/3} + j^{1/3})^2|i^{2/3} - j^{2/3}|\text{erf}\Big(\dfrac{|i^{2/3} - j^{2/3}|}{\sqrt{i^{1/3} + j^{1/3}}} \Big)$ & $7 (0.02)$& $14 (0.02)$ & $4 (0.03)$& $ 8 (0.05)$ & modified flux reaction rates \cite{osinsky2024hydrodynamic}\\
    \hline

    \rule{0pt}{1.5em} $(i^{1/3} + j^{1/3})^{5/2}\exp\Big[{-\dfrac{(i^{2/3} - j^{2/3})^2}{i^{1/3} + j^{1/3}}}\Big]$ & $7 (0.02)$ & $13 (0.02)$ & $5 (0.03)$ & $8 (0.05)$ & modified flux reaction rates \cite{osinsky2024hydrodynamic}\\
    \hline
    
    \rule{0pt}{1.7em}$(i^{2/3} + j^{2/3})\sqrt{(i^{2/9} + j^{2/9})}\exp\Big[-c\Big(\dfrac{i^{1/3}j^{1/3}}{i^{1/3} + j^{1/3}}\Big)^{4}\Big]$& $5(0.01)$ & $10(0.01)$& $4 (0.03)$& $8 (0.05)$& emulsion coalescence\cite{coulaloglou1977description, koch2008h} \\
    \hline
    \rule{0pt}{1.7em}$(i^{1/3} + j^{1/3})^2\sqrt{\dfrac{1}{i} + \dfrac{1}{j}}$& $4 (0.01)$& $7(0.01)$ & $3 (0.03)$ & $6 (0.04)$ & ballistic kernel \cite{Leyvraz03} \\
    \hline

    \rule{0pt}{1.7em}$K_{ij}=(i^{1/3} + j^{1/3})^2\cdot 
    \left|\dfrac{1}{i} - \dfrac{1}{j} \right|$& $3 (0.01)$ & $6 (0.01)$ & $3 (0.03)$ & $5 (0.04)$& modified ballistic kernel \cite{zagidullin2022aggregation} \\
    \hline

    \rule{0pt}{1.7em} $(i^{1/3} + j^{1/3})^2\cdot \left|\dfrac{1}{1 - \mathbf{i}\cdot i^{2/3}} - \dfrac{1}{1 - \mathbf{i}\cdot j^{2/3}}\right|$ & $4 (0.01)$ & $4 (0.01)$ & $3 (0.03)$& $4 (0.04)$ & orthokinetic interaction \cite{sheng2007simulation, khmelev2024smoke, khmelev2024mathematical} \\
    \hline

    \rule{0pt}{1.7em} $\dfrac{i^{2/3}j^{2/3}}{i^{1/3} + j^{1/3}}$ & $3 (0.01) $ & $6 (0.01)$& $3 (0.03)$ & $5 (0.04) $& hydrodynamic interaction \cite{sheng2007simulation, khmelev2024mathematical} \\
    \hline

    \rule{0pt}{1.7em} $(i^{1/3} + j^{1/3})^2\sqrt{i^{1/3}~j^{1/3}}\exp\Big[-c(i^{1/3}j^{1/3})^{1/3}\Big]$ & $4 (0.01)$ & $6(0.01)$ & $3 (0.03)$ & $6 (0.04)$ & coalescence of ﬂuid particles
(bubbles and drops) \cite{liao2010literature, konno1988coalescence} \\
    \hline

    \rule{0pt}{1.7em} $(i^{1/3} + j^{1/3})^2\sqrt{(i^{2/3} + j^{2/3})}$ & $3 (0.01)$ & $6 (0.01)$& $3 (0.03)$ & $6 (0.04)$ & coalescence of ﬂuid particles
(bubbles and drops) \cite{luo1995coalescence} \\
    \hline
\end{tabular}
\caption{Maximal rank of a block (compression rate \%) of different kernels in the mosaic--skeleton format. The constant $c$ can have different values in different kernels. For more details, see the source in the list of references.}
\label{table:mosaic_ranks}
\end{table}

%% file: calculation_operators.tex
\section{Mosaic--skeleton acceleration of Smoluchowski operator computation}\label{sec:smol_operations}


A common way to solve a system of ODE numerically is to integrate it with some finite difference scheme. In particular, a family of Runge--Kutta methods allows to utilize the adaptive time-steps \cite{matveev2024adaptive}. 
In our case, the most time consuming operation is the evaluation of the right--hand side of the equation. Hence, this Section is devoted to its reduction.

The Smoluchowski operator consists of two parts: $f_1(s)$ and $f_2(s)$ for $s = 1, \ldots, M$, defined in Eq.~\eqref{eq:smoluch_euler_approximation}.
In Section~\ref{sec:lowrank_kernels} we show how the low--rank structure of a kernel can be used for their fast evaluation. Such trick reduces the complexity of evaluation of $f_1(s)$ and $f_2(s)$ from quadratic to sublinear.
Further, we demonstrate how to accelerate similar computations but for kernels allowing low-parametric representation in mosaic--skeleton format, described in the previous Section~\ref{sec:mosaic_skeleton_kernels}. Utilization of this structure allows to get almost the same asymptotic acceleration but for the broader family of kernels.

For convenience, we use the notations in agreement with the original paper \cite{tyrtyshnikov1996mosaicskeleton}.
Let the mosaic matrix $K$ be composed of blocks $B_p$, $p = 1, \ldots, P$.
Denote $\Gamma(B_p)$ the matrix that coincides with $K$ on the block $B_p$ and has zeroes elsewhere.
Then we get the following representation for the kernel:
\begin{equation} \label{eq:kernel_mosaic_sum}
    K = \sum\limits_{p = 1}^{P} \Gamma(B_p).
\end{equation}
Next, for the mosaic--skeleton matrix $K$ we define its mosaic rank:
\begin{equation} \label{eq:def_mosaic_rank}
    \mrank K = \dfrac{1}{2M} \mem K = \dfrac{1}{2M} \sum\limits_{p = 1}^{P} \mem B_p,
\end{equation}
where $\mem B_p$ is the number of memory cells required for storing the representation of the block $B_p$. Here we assume that all blocks are square.
Thus, for the block of size $m_p \times m_p$ and $\rank B_p = r_p$ its number of parameters equals $\mem B_p = \min(m_p^2, 2 r_p m_p)$.
This concept extends the classical definition of the matrix rank in terms of the number of parameters stored and also matrix-by-vector multiplication complexity, see Section~\ref{sec:smol_operation:matvec}.

As soon as definition is given, we need to find a relation between the mosaic rank and the matrix size $M$. The common approach is to assume that the rank of every block (except dense ones) is bounded by some constant:
\begin{equation*}
    \rank B_p \leq R.
\end{equation*}
For our mosaic partitioning types (see Figure~\ref{fig:different_splitting}) the total number of blocks of size $\dfrac{M}{2^k} \times \dfrac{M}{2^k}$ is proportional to $2^k$ for every $k = 1, \ldots, L$.
Here $L = \log_2 M - k_0$ is the number of rows and columns bisection levels.
Then the low--rank blocks require not more than
\begin{equation*}
    \sum\limits_{k=1}^{L} 2^k \cdot R \dfrac{M}{2^k} = \mathcal{O}(R M \log M)
\end{equation*}
memory cells, where $R$ is the maximal rank among all blocks.
The number of dense blocks of size $\dfrac{M}{2^{L}} \times \dfrac{M}{2^{L}}$ is also proportional to $2^{L}$ (they can be of the smallest size only), which requires the following number of parameters for them:
\begin{equation*}
    \mathcal{O}\left(2^L \cdot \left(\dfrac{M}{2^L}\right)^2\right) = \mathcal{O}(2^{k_0} M) = \mathcal{O}(M).
\end{equation*}

Thus, the total number of parameters for the mosaic--skeleton representation of the kernel ($\mem K$) is bounded by $\mathcal{O}(R M \log M)$, and from \eqref{eq:def_mosaic_rank} it follows that the mosaic rank of $K$ is small:
\begin{equation} \label{eq:mosaic_rank_bound}
    \mrank K = \mathcal{O}(R \log M).
\end{equation}

For a vector $v$ we also denote its slice of size $l$, starting from the index $\alpha$:
\begin{equation*}
    v(\alpha : \alpha + l) \equiv [v_\alpha, v_{\alpha + 1}, \ldots, v_{\alpha + l - 1}]^\top.
\end{equation*}
Further, matrix and vector indexing starts from one.

\subsection{Second operator, matrix--vector multiplication} \label{sec:smol_operation:matvec}

At first, we discuss the second component of the Smoluchowski operator because its computation is quite simple:
\begin{equation*}
     f_2(s) = n_s \sum\limits_{j = 1}^{M} K_{sj} n_j.
\end{equation*}
In the matrix--vector form, we need to evaluate:
\begin{equation*}
    f_2 = n \odot (K n),
\end{equation*}
where $\odot$ denotes element-wise product.

Using representation~\eqref{eq:kernel_mosaic_sum}, we get the following relation:
\begin{equation*}
    f_2 = n \odot \left(\sum\limits_{p = 1}^{P} \Gamma(B_p)\right) n = n \odot \left(\sum\limits_{p = 1}^{P} \Gamma(B_p) n\right).
\end{equation*}
In fact, calculating every summand $\Gamma(B_p) n$, $p = 1, \ldots, P$, results in the following local matrix-vector product and addition:
\begin{equation*} \label{eq:matvec_mosaic}
    f_2(\alpha_p : \alpha_p + m_p) \gets f_2(\alpha_p : \alpha_p + m_p) + B_p n(\beta_p : \beta_p + m_p).
\end{equation*}
Or in index form:
\begin{equation*}
    f_2(\alpha_p + s - 1) \gets f_2(\alpha_p + s - 1) + \sum\limits_{i=1}^{m_p} B_p(s, i) n(\beta_p + i - 1), \quad s = 1, \ldots, m_p.
\end{equation*}
Here $(\alpha_p, \beta_p)$ is the upper left corner of a block $B_p$ of size $m_p \times m_p$.
Then getting
\begin{equation*}
    f_2 \gets f_2\odot n
\end{equation*}
finalizes the evaluation of the operator $f_2$.

For a dense block, its direct multiplication by a vector costs $m_p^2$ operations.
And for the low--rank one, it costs $2r_p m_p$ operations if FMA (Fast Multiply plus Add) is used, as was discussed in Section~\ref{sec:lowrank_kernels}.
So, for both of the types of block $B_p$ this procedure requires exactly $\mem B_p$ operations.
Summation over all $P$ blocks leads to the total complexity equal to $\mem K$.
Combining \eqref{eq:def_mosaic_rank} and \eqref{eq:mosaic_rank_bound}, we obtain:
\begin{equation} \label{eq:complexity_matvec}
    \mathcal{C}_{f_2}(M) = \mathcal{O}(R M \log M),
\end{equation}
where $\mathcal{C}_{f_2}(M)$ is the number of arithmetic operations required for computation of the operator $f_2$ over a vector of size $M$ using the mosaic--skeleton format.
Recall that for a low--rank kernel we had $\mathcal{O}(RM)$.
The new complexity is only on a logarithmic factor higher, while it allows to consider a broader class of kernels.

\subsection{First operator, convolution}


The first term in Eq.~\eqref{eq:smoluch_euler_approximation} is defined as follows:
\begin{equation*}
    f_1(s) = \sum\limits_{i + j = s} K_{ij} n_i n_{j}, \quad s = 2, \ldots, M.
\end{equation*}
Let us denote $z = \conv(y, A, x)$ the full convolution of two vectors with a kernel:
\begin{equation} \label{eq:conv_direct}
    z_{s - 1} = \sum\limits_{i + j = s} y_i K_{ij} x_j, \quad s = 2, \ldots, 2M.
\end{equation}
Here the resulting vector~$z$ has the size $2M - 1$.
For obtaining $f_1$ we need to evaluate its first $M-1$ components: $f_1 = \conv(n, K, n)(1 : M)$.

For the kernel $K$ in the mosaic--skeleton format~\eqref{eq:kernel_mosaic_sum}, convolution transforms into the following relation:
\begin{equation*}
    f_1(s) = \sum_{p = 1}^{P} \sum\limits_{i + j = s} \Gamma(B_p)_{ij} n_i n_{j}.
\end{equation*}
As far as $\Gamma(B_p)$ has zeroes elsewhere except for positions $(\alpha_p : \alpha_p + m_p) \times (\beta_p : \beta_p + m_p)$ of the block $B_p$, each of $P$ summands can be calculated using the following formula:
\begin{equation*}
    f_1(\gamma_p : \gamma_p + l_p) \gets f_1(\gamma_p : \gamma_p + l_p) + \conv\left( n(\alpha_p : \alpha_p + m_p), B_p, n(\beta_p : \beta_p + m_p) \right),
\end{equation*}
where $\gamma_p = \alpha_p + \beta_p$ and $l_p = 2 m_p - 1$.
If $\gamma_p + l_p > M + 1$ then the addition is performed till the index $M$, or omitted at all for $\gamma_p > M$.

Convolutions with the dense block $B_p$ is performed as it is defined in Eq.~\eqref{eq:conv_direct} requiring $2 m_p^2 = 2 \mem B_p$ operations.
For the low--rank blocks we utilize the FFT-based procedure in the same way as it is described in Section~\ref{sec:lowrank_kernels} above.
It requires $\mathcal{O}(r_p m_p \log m_p) = \mathcal{O}(\mem B_p \log m_p) = \mathcal{O}(\mem B_p \log M)$ operations if $\rank B_p = r_p$.
Summation over all $P$ dense and low--rank blocks leads us to $\mathcal{O}(\mem K \log M)$ complexity bound.
Using Eqs.~\eqref{eq:def_mosaic_rank}~and~\eqref{eq:mosaic_rank_bound}, we obtain the total complexity of operator $f_1$ evaluation with a kernel in the mosaic--skeleton format:
\begin{equation} \label{eq:complexity_conv}
    \mathcal{C}_{f_1}(M) = \mathcal{O}(R M \log^2 M),
\end{equation}
the number of arithmetic operations for vector of size $M$.
Recall that in the case of a low--rank kernel we had $\mathcal{O}(R M \log M)$, which is also now being increased by the same logarithmic factor as $f_1$.

Thus, combining Eqs.~\eqref{eq:complexity_matvec}~and~\eqref{eq:complexity_conv} we obtain the total complexity for the Smoluchowski operator:
\begin{equation} \label{eq:complexity_total}
    \mathcal{C}_{f_1 + f_2}(M) = \mathcal{O}(R M \log^2 M),
\end{equation}
for a kernel $K$ in the mosaic--skeleton format with off-diagonal blocks rank bounded by $R$.
There is a wide range of full rank kernels (see Table~\ref{table:mosaic_ranks}) used in practice, for which we can gain significant speedup while solving the Smoluchowski equations.

%% file: numerical_experiments.tex
\section{Numerical experiments}

In this section we present the results of several numerical experiments allowing to investigate the properties of our algorithm.
Its convergence rate and computational complexity are studied.
The experiments have been carried out for various scenarios with two types of the initial particle size distributions: the monodisperse and distributions with a wide spectrum of sizes.

Further, we denote \texttt{FDMSk} our finite difference scheme with the mosaic-skeleton approximation of a kernel.
The label \texttt{LRMC} corresponds to the efficient majorant Monte Carlo method~\cite{osinsky2024low} that we use for comparison.
The alternative method also utilizes the low-rank representations of the kernel coefficients in terms of the the special majorant functions. 
It shows a good performance relying on the rank $R$ of majorant function: $\mathcal{O}\left(R\log M\right)$ per collision.
For more details and implementation techniques we refer to the original paper~\cite{osinsky2024low}.

In our experiments, we consider three kernels (constant, \eqref{eq:kernel_stream} and \eqref{eq:kernel_baikal}):
\begin{align*}
    K_{i j}^1 &= 2, &
    K_{i j}^2 &= (i^{1/3} + j^{1/3})^2|i^{2/3}-j^{2/3}|, &
    K_{i j}^3 &= \frac{(i+j)\left(i^{1 / 3}+j^{1 / 3}\right)^{2 / 3}}{(i j)^{5 / 9}\left|i^{2 / 3}-j^{2 / 3}\right|}.
\end{align*}

We obtain the low-rank majorant functions for all these kernels.
The first one has rank $R = 1$ itself, the second has rank $R = 2$ and the third one can be estimated by rank $R = 6$ function:
\begin{align*}
    K_{ij}^2 &\le \hat{K}_{ij}^2 = 2 (i^{4/3}+j^{4/3}), & R &= 2, \\
    ~ & ~\\
    K_{ij}^3 &= \dfrac{(i^{4/3} + i^{2/3}j^{2/3} + j^{4/3})(i^{1/3} + j^{1/3})^{2/3}}{(ij)^{5/9}|i-j|} \le \hat{K}_{ij}^3 = \dfrac{(i^{4/3} + i^{2/3}j^{2/3} + j^{4/3})(i^{2/9} + j^{2/9})}{(ij)^{5/9}}, & R &= 6.
\end{align*}

In all experiments, we approximate blocks of the kernel $K_{ij}$ using the mosaic partitioning from the Figure~\ref{fig:different_splitting:tridiag} and matrix cross method~\cite{zheltkov2015parallel} with the relative accuracy $10^{-6}$. 
The values of the maximal ranks of blocks are given in the Table~\ref{table:mosaic_ranks}.
In our benchmarks, we add the time of construction of the approximation for the kernel to the total running time of the algorithm.

We perform our numerical experiments on a single node of the cluster of INM~RAS~\cite{INM_cluster}.
All calculations for the \texttt{FDMSk} were performed on a single thread.
The running time of the \texttt{LRMC} algorithm in the Table~\ref{table:convergenceanalysis1} recalculated, as if all calculations were performed on the single thread, due to the fact that we used averaging of several solutions to increase accuracy
The implementation of the \texttt{FDMSk} algorithm can be found by the link \footnote[1]{\href{https://github.com/DrEternity/FDMSk-Smoluchowski}{https://github.com/DrEternity/FDMSk-Smoluchowski}}.

\subsection{Basic benchmark}

Unfortunately, there are only a few know solutions in the analytical form for the Smoluchowski equations.
For example, for the monodisperse initial conditions $n_k(t=0) = \delta_{k,1}$ with the kernel $K_{ij}^1=2$ such solution has been found by Smoluchowski \cite{Smol17}:
\begin{equation} \label{eq:theor_solution}
    n_k(t)=\frac{1}{(1+t)^2}e^{-(k-1)\ln(1+1/t)}.
\end{equation}
We verify that our algorithm allows to obtain an accurate numerical solution of this problem at the time point $t = 100$. We present these results in Table~\ref{table:const_kernel_results} demonstrating the dependency of the solution error on the number of equations and the time step size.
The numerical error is measured in terms of the first moment $M_1$, defined the following way:
\begin{equation} \label{eq:m1-norm}
    \|n^{\texttt{FDMSk}} - n^{\texttt{Theory}}\|_{M_1} = \sum_{i} i \cdot |n^{\texttt{FDMSk}}_i - n^{\texttt{Theory}}_i|.
\end{equation}
The accuracy increases with both of the parameters refinement.
One may note that for some cases the error stagnates.
That is because it has two components, corresponding to the mass leaking and finite difference scheme error.
They are controlled by the maximal particle size and the time step, respectively.
Each of them gives its own limitations, so we need to refine them simultaneously.

However, the advantages of our method cannot be fully demonstrated in such simple an example: a lot of methods show great performance basing on this classical benchmark but lose their efficiency in cases with more complicated kinetic coefficients. 
Below, we consider more complex kernels without known analytical solutions.

\begin{table}[ht!]
\centering
\begin{tabular}[c]{| c || c | c | c | c |}
     \hline
     \diagbox{$M$}{$\Delta t$} & $0.5$ & $0.1$ & $0.05$ & $0.01$\\    
     \hline
    $256$ & $3 \cdot 10^{-2}$ & $3 \cdot 10^{-2}$ & $3 \cdot 10^{-2}$ & $3 \cdot 10^{-2}$ \\
    $512$ & $2 \cdot 10^{-3}$ & $2 \cdot 10^{-3}$ & $2 \cdot 10^{-3}$ & $2 \cdot 10^{-3}$\\
    $1024$ & $2 \cdot 10^{-5}$ & $6 \cdot 10^{-6}$ & $6 \cdot 10^{-6}$ & $6 \cdot 10^{-6}$\\
    $4096$ & $1 \cdot 10^{-5}$ & $2 \cdot 10^{-7}$ & $2 \cdot 10^{-7}$ & $2 \cdot 10^{-7}$\\
    $16384$ & $9 \cdot 10^{-6}$ & $1 \cdot 10^{-8}$ & $9 \cdot 10^{-9}$ & $9 \cdot 10^{-9}$\\
    $65536$ & $9 \cdot 10^{-6}$ & $2 \cdot 10^{-9}$ & $1 \cdot 10^{-9}$ & $9 \cdot 10^{-10}$\\

     \hline
    
    \hline
\end{tabular}
\caption{Relative error of \texttt{FDMSk} method in norm $\|\cdot\|_{M_1}$  for case (\ref{eq:theor_solution}) and $t=100$ with $M$ equations and time step $\Delta t$.}
\label{table:const_kernel_results}
\end{table}

\subsection{Verification with Monte Carlo}

In this subsection, kernels $K_{ij}^2, K_{ij}^3$ are studied. 
We have already noted that they have a more complex structure in comparison with the constant kernel.
Since there is no available analytical solution for both of these kernels, we investigate the convergence of our algorithm numerically fixing $t = 100$ for $K_{ij}^3$ and $t = 0.5$ in case of $K_{ij}^2$. 
In order to do this, we check the mutual convergence of the methods of different nature: \texttt{FDMSk} and \texttt{LRMC}.
The value $\varepsilon_1$ denotes the mutual error, measured in terms of the first moment (see Eq.~\eqref{eq:m1-norm}):
\begin{equation*}
    \varepsilon_1 = \dfrac{\left\|n^{\texttt{FDMSk}} - n^{\texttt{LRMC}}\right\|_{M_1}}{\left\|n^{\texttt{Ref}}\right\|_{M_1}} = \dfrac{\sum_{i} i \cdot |n_{i}^{\texttt{FDMSk}} - n_{i}^{\texttt{LRMC}}|}{1}.
\end{equation*}
Here, the first moment of solution is considered unitary since there is no leak of mass.
To evaluate convergence for a particular method, we also use the second moment integral characteristic:
\begin{equation*}
    \|n\|_{M_2} = \sum\limits_{i} i^2 n_i.
\end{equation*}
Thus, the convergence is captured by the second moment of the reference solution: 
\begin{equation*}
    \varepsilon_2 = \dfrac{\left| \left\|n^{\texttt{Numerical}}\right\|_{M_2} - \left\|n^{\texttt{Ref}}\right\|_{M_2}\right|}{\left\|n^{\texttt{Ref}}\right\|_{M_2}}.
\end{equation*}
The accuracy of the numerical solution obtained by the different methods depends on different parameters.
For \texttt{FDMSk} the accuracy depends on the number of equations $M$ and the step of the difference scheme $\Delta t$. 
And for the \texttt{LRMC} it depends on the number of particles $V$ supported in the system and the number of simulations $L$ used for averaging~\cite{goodson2002efficient}.
In order to obtain the reference solutions we take $M = 2^{12}$, and $\Delta t^2 =10^{-5}$, $\Delta t^3 = 10^{-2}$ for $K_{ij}^2$, $K_{ij}^3$ respectively in case of \texttt{FDMSk}.
For the \texttt{LRMC} we setup reference parameters $V=10^7$ and $L=10^6$.

\begin{figure}[ht]
    \centering
    \includegraphics[width=0.8\textwidth]{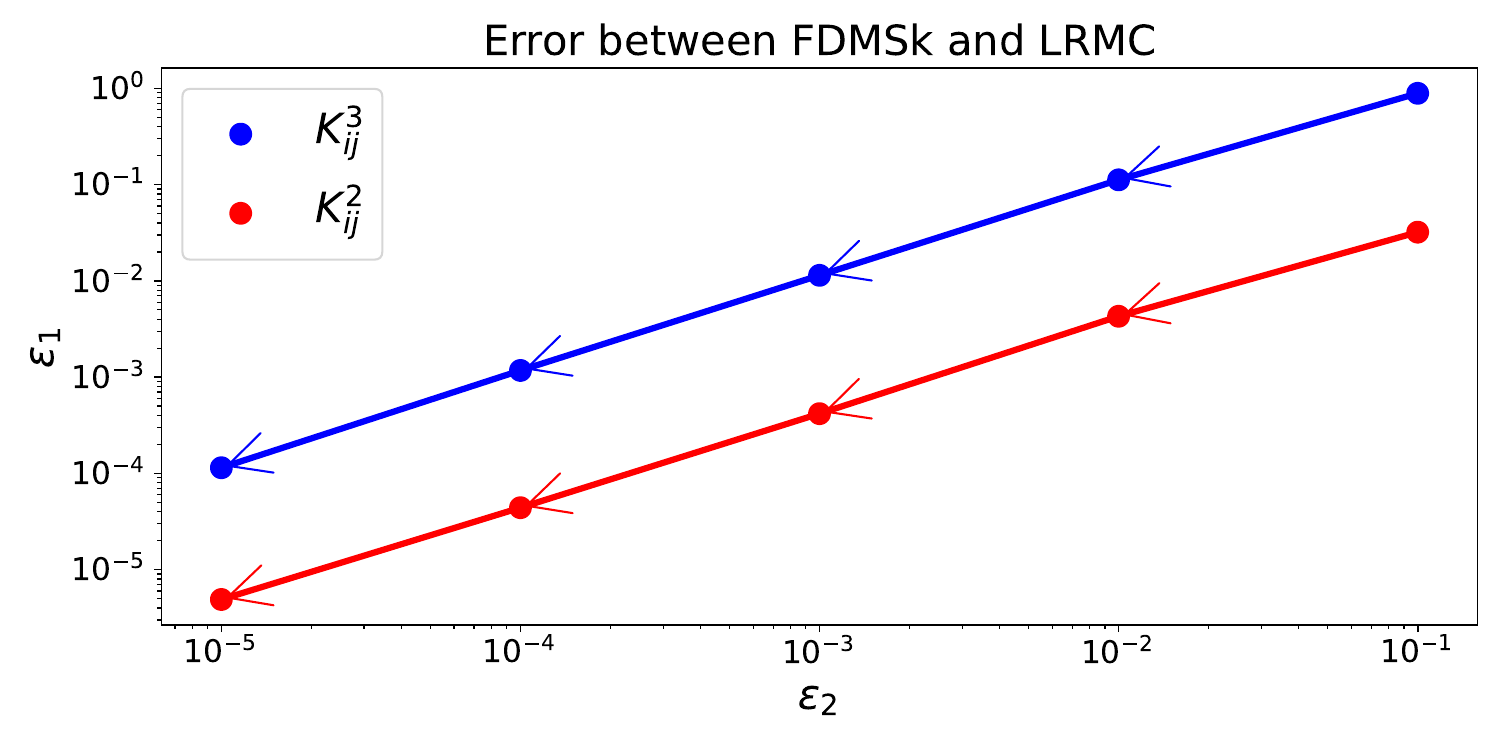}
    \caption{The mutual error $\varepsilon_1$ of the \texttt{FDMSk} and \texttt{LRMC} with increasing accuracy $\varepsilon_2$ of each of them for kernels $K_{ij}^{2}, K_{ij}^3$. Linear decreasing of mutual error is observed.}
    \label{fig:mutual_error}
\end{figure}

In this numerical experiment, we again consider a monodisperse initial condition $n_k(t=0) = \delta_{k,1}$.
For both \texttt{FDMSk} and \texttt{LRMC} we setup such simulation parameters that allow to achieve the required convergence rates~$\varepsilon_2$ and measure the mutual error~$\varepsilon_1$ of the obtained solutions.
From the Figure~\ref{fig:mutual_error}, we see that as $\varepsilon_2$ decreases, both algorithms demonstrate the convergence to almost the same solution. 
The value $\varepsilon_1$ tends to zero with the same rate.
This convinces us of the correctness of the \texttt{FDMSk} technique.

Here, we also note the advantages of the finite difference scheme over Monte Carlo methods (see Table~\ref{table:convergenceanalysis1}).
If the high accuracy of solution is necessary (e.g., better than $10^{-2}$), the novel \texttt{FDMSk} approach achieves it significantly faster than \texttt{LRMC}.

\begin{table}[ht!]
    \begin{center}
         \begin{tabular}[c]{| c | c | c || c | c |}
         \hline
        \multirow{2}{*}{$\varepsilon_2$} & \multicolumn{2}{c||}{$K_{ij}^2$}& \multicolumn{2}{c|}{$K_{ij}^3$}\\
    
         & LRMC [sec] & FDMSk [sec] & LRMC [sec] & FDMSk [sec] \\
         \hline 
         \hline
        $10^{-1}$ & 0.28 & 0.7 & 0.034 & 9.76\\
        $10^{-2}$ & 9.2 & 2.1 & 3.6 & 10.569\\
        $10^{-3}$ & 821 & 7.5 & 315.6 & 23.821\\
        $10^{-4}$ & 77~975 ($\approx$ 21 hours) & 8.2 & 27 192 ($\approx$ 7.5 hours) & 54.490\\
        $10^{-5}$ & 7~908~827 ($\approx$ 91 days) & 30.08 & 2 599 938 ($\approx$ 30 days) & 129.692\\
        \hline 
    	\end{tabular}
    \end{center}
    \caption{The running time of the algorithms under different accuracy requirements $\varepsilon_2$}
    \label{table:convergenceanalysis1}
\end{table}


The other important advantage of \texttt{FDMSk} is its efficiency while operating with a large number of equations.
It becomes necessary if we have a wide range of particle sizes. 
For example, consider a system with the following initial conditions:
\begin{equation*}
    n_k(t=0) = \left\lbrace
    \begin{matrix}
        & \dfrac{1}{k + 1}, & k \leq M, \\
        & ~ & ~ \\
        & 0, & k > M. 
    \end{matrix}\right.
\end{equation*}
Here, or testing on a computationally intensive case we set the maximal particle size as $M = 10^6$ and consider the target dimensionless time stamp $t = 1$. 

The disadvantage of the \texttt{LRMC} method is its ambiguity in choosing a low-rank majorant for the kernel. For many well-known kernels, it is easy to write out a majorant that will asymptotically allow a high acceptance rate. However, for example, for the $K_{ij}^3$, we hypothesize that it is impossible to construct an asymptotically correct low-rank majorant. Unfortunately, our majorant $\hat{K}_{ip}^3$ has a poor acceptance rate:
\begin{equation*}
\lim_{i \rightarrow +\infty} P(Accept \; | \;\{i, p\}) = \lim_{i \rightarrow +\infty}\frac{K_{ip}^3}{\hat{K}_{ip}^3} = \lim_{i \rightarrow +\infty} \dfrac{1}{i} = 0, \quad p \in \mathbb{N}
\end{equation*}

This fact is confirmed numerically by the experiment.
The running time of the \texttt{FDMSk} algorithm with adaptive time-step selection in the Runge-Kutta scheme takes \textbf{2782 seconds}. 
Meanwhile, one run of the \texttt{LRMC} with $V = 10^6$ simulation requires \textbf{53 hours} and the acceptance rate is $2\cdot10^{-5}$. 
This demonstrates the best side of the proposed \texttt{FDMSk} method and demonstrates its extremely good properties in solving truly large systems.

\subsection{Complexity analysis}

In this section we verify the theoretical estimates of the complexity of the \texttt{FDMSk} algorithm in practice. 
Recall that the estimate~\eqref{eq:complexity_total} on the algorithmic complexity for each time step of the \texttt{FDMSk} method is
\begin{equation*}
    \mathcal{C}_{f_1 + f_2}(M) = \mathcal{O}\left(RM\log^2M\right).
\end{equation*}
It consists of evaluation of the operators $f_1$ and $f_2$ with complexities $\mathcal{O}(RM \log^2 M)$ and $\mathcal{O}(RM \log M)$, respectively, see Eqs.~\eqref{eq:complexity_conv} and~\eqref{eq:complexity_matvec}.

\begin{figure}[ht]
    \centering
    \begin{subfigure}{0.45\textwidth}
      \centering
      \includegraphics[width=1.0\linewidth]{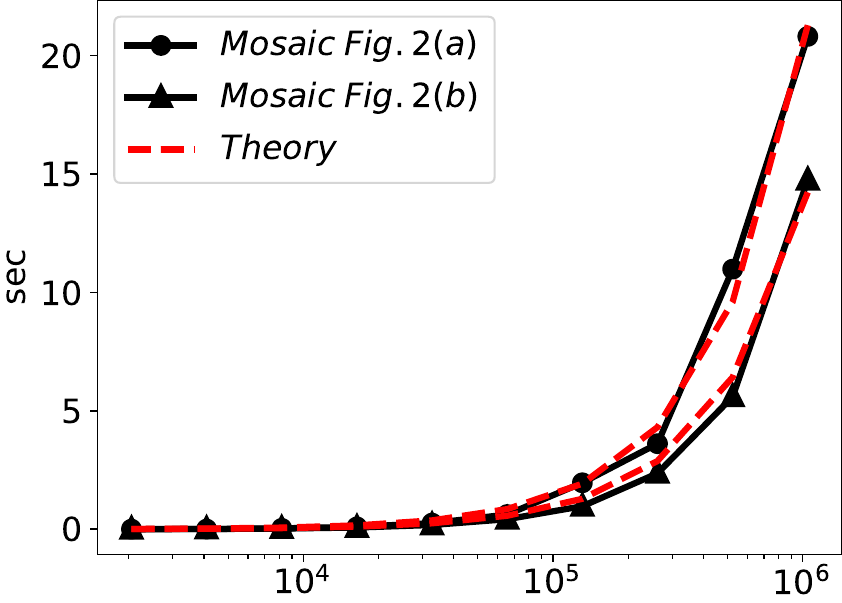}\\
      \caption{$f_1$}
    \end{subfigure}%
    \hfill
    \begin{subfigure}{0.47\textwidth}
      \centering
      \includegraphics[width=1.0\linewidth]{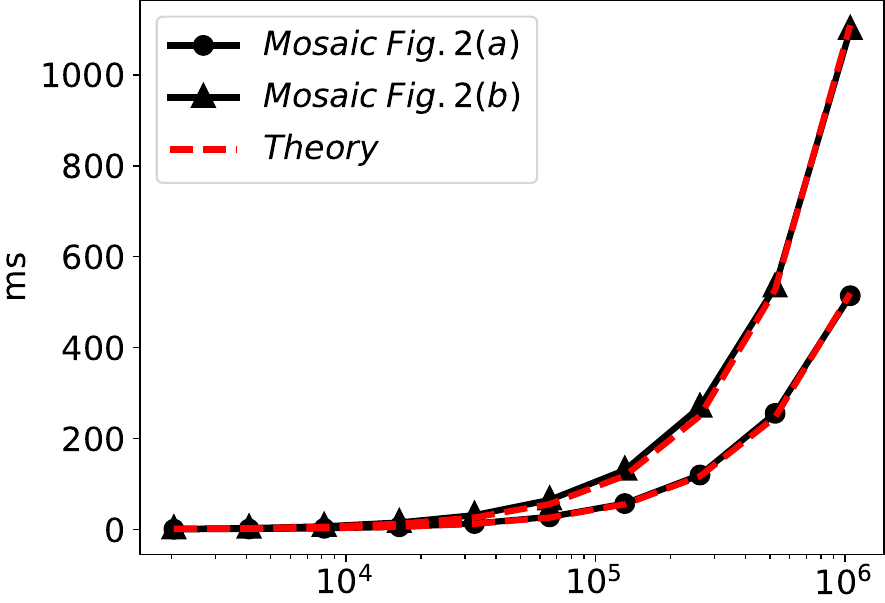}\\
      \caption{$f_2$}
    \end{subfigure}
    \caption{Evaluation time of the operators $f_1$ and $f_2$ in mosaic-skeleton format compared to the theoretical estimates.}
    \label{fig:check_complexity}
\end{figure}

In the experiment we study the performance of our code on the evaluation of $f_1$ and $f_2$ in the case of kernel $K^{3}_{ij}$. On Figure~\ref{fig:check_complexity} an excellent agreement of the theoretical estimates with the actual performance of the algorithm is demonstrated for two different mosaic partitionings.
Precise running time values can be found in Table~\ref{table:check_complexity}. 
Note that calculations with $M = 2^{20}$ require modest amount of time and can be performed on a personal computer.

\begin{table}[ht]
    \begin{center}
         \begin{tabular}[c]{| c | c | c | c | c |}
         \hline
         M & $f_1$ [sec] Figure~\ref{fig:different_splitting:monodiag} & $f_2$ [ms] Figure~\ref{fig:different_splitting:monodiag} & $f_1$ [sec] Figure~\ref{fig:different_splitting:tridiag} & $f_2$ [ms] Figure~\ref{fig:different_splitting:tridiag} \\
         \hline
         \hline
         $2^{11}$ & 0.0048 & 0.3580 & 0.0022 & 0.8576 \\
         $2^{12}$ & 0.0078 & 0.8352 & 0.0081 & 2.7439 \\
         $2^{13}$ & 0.0377 & 2.2932 & 0.0290 & 6.5832 \\
         $2^{14}$ & 0.1002 & 5.6186 & 0.0728 & 15.080 \\
         $2^{15}$ & 0.2591 & 12.820 & 0.1796 & 31.206 \\
         $2^{16}$ & 0.6207 & 27.435 & 0.4245 & 64.712 \\
         $2^{17}$ & 1.9605 & 56.718 & 0.9710 & 130.33 \\
         $2^{18}$ & 3.6114 & 119.53 & 2.3615 & 269.00 \\
         $2^{19}$ & 10.980 & 255.25 & 5.5937 & 533.12 \\
         $2^{20}$ & 20.804 & 514.06 & 14.737 & 1099.9 \\
        \hline 
    	\end{tabular}
    \end{center}
    \caption{Operation time in mosaic-skeleton format for the operators $f_1$, $f_2$ for $K_{ij}^3$.}
    \label{table:check_complexity}
\end{table}